\documentclass{amsart}
\usepackage{amssymb,amsfonts}
\usepackage{enumerate}
\usepackage{amsmath,calligra,mathrsfs}
\usepackage{graphicx}
\usepackage{mathtools}
\usepackage{leftidx}
\usepackage{tikz-cd}
\usepackage{bm}
\usepackage{url}
\usepackage{makecell}
\usepackage{underoverlap}

\newtheorem{thm}{Theorem}[section]
\newtheorem{cor}[thm]{Corollary}
\newtheorem{prop}[thm]{Proposition}
\newtheorem{lem}[thm]{Lemma}

\theoremstyle{definition}
\newtheorem{defn}[thm]{Definition}

\newtheorem{exmp}[thm]{Example}

\newtheorem{chkpt}[thm]{Checkpoint}

\theoremstyle{remark}
\newtheorem{rem}[thm]{Remark}

\makeatletter
\let\c@equation\c@thm
\makeatother
\numberwithin{equation}{section}

\def\bthm{\begin{thm}}
\def\ethm{\end{thm}}
\def\blm{\begin{lem}}
\def\elm{\end{lem}}
\def\bdf{\begin{defn}}
\def\edf{\end{defn}}
\def\bpf{\begin{proof}}
\def\epf{\end{proof}}
\def\bpp{\begin{prop}}
\def\epp{\end{prop}}
\def\bcor{\begin{cor}}
\def\ecor{\end{cor}}
\def\brm{\begin{rem}}
\def\erm{\end{rem}}
\def\beg{\begin{exmp}}
\def\eeg{\end{exmp}}
\def\bD{\mathbb{D}}

\def\bG{\mathbb{G}}

\def\bN{\mathbb{N}}

\def\bP{\mathbb{P}}
\def\bQ{\mathbb{Q}}

\def\bZ{\mathbb{Z}}
\def\cA{\mathcal{A}}

\def\cC{\mathcal{C}}
\def\cD{\mathcal{D}}
\def\cE{\mathcal{E}}
\def\cF{\mathcal{F}}

\def\cH{\mathcal{H}}

\def\cL{\mathcal{L}}

\def\cO{\mathcal{O}}

\def\scL{\mathscr{L}}

\def\frX{\mathfrak{X}}

\newcommand{\raq}{\,\rightarrow \,}

\newcommand{\xraq}[2][]{\, \xrightarrow[#1]{#2} \,}
\newcommand{\xlaq}[2][]{\, \xleftarrow[#1]{#2} \,}

\newcommand{\ra}{\rightarrow}

\newcommand{\rinto}{\hookrightarrow}

\newcommand{\xra}[2][]{\xrightarrow[#1]{#2}}
\newcommand{\xla}[2][]{\xleftarrow[#1]{#2}}

\newcommand{\Mod}{{\rm Mod}}

\newcommand{\Ab}{{\rm Ab}}

\newcommand{\Tr}{{\rm Tr}}

\newcommand{\Tor}{{\rm Tor}}

\newcommand{\RHom}{{\bm R}{\rm Hom}}
\newcommand{\RHomcom}{{\bm R} \underline{{\rm Hom}}}

\newcommand{\Ob}{\rm{Ob}}
\newcommand{\op}{{\rm op}}

\newcommand{\Homcom}{\underline{{\rm Hom}}}

\newcommand{\cHom}{\mathscr{H}\text{\kern -3pt {\calligra\large om}}\,}
\newcommand{\RcHom}{{\bm R}\cHom}

\newcommand{\Spec}{{\rm Spec}}
\newcommand{\Proj}{{\rm Proj}}
\newcommand{\uProj}{\underline{{\rm Proj}}}

\newcommand{\QCoh}{{\rm QCoh}}
\newcommand{\Coh}{{\rm Coh}}

\newcommand{\Pic}{{\rm Pic}}

\newcommand{\Dbcoh}{\cD^b_{{\rm coh}}}

\newcommand{\Dqcoh}{\cD_{{\rm qcoh}}}

\newcommand{\Dmcoh}{\cD^{-}_{{\rm coh}}}

\newcommand{\Dperf}{\cD_{{\rm perf}}}

\newcommand{\perf}{{\rm perf}}

\newcommand{\Dpcoh}{\cD^{+}_{{\rm coh}}}

\newcommand{\Ex}{{\rm Ex}}
\newcommand{\WDiv}{{\rm WDiv}}

\newcommand{\GrcA}{{\rm Gr} ( \cA)}
\newcommand{\GrA}{{\rm Gr} (A)}
\newcommand{\Gr}{{\rm Gr}}
\newcommand{\Torp}{{\rm Tor}^+}
\newcommand{\Iptr}{I^+\text{-triv}}
\newcommand{\Imtr}{I^-\text{-triv}}
\newcommand{\RGIp}{{\bm R}\Gamma_{I^+}}
\newcommand{\CIp}{\check{\cC}_{I^+}}
\newcommand{\RGIm}{{\bm R}\Gamma_{I^-}}
\newcommand{\CIm}{\check{\cC}_{I^-}}

\DeclareMathAlphabet{\mathpzc}{OT1}{pzc}{m}{it}

\title{Survey on homological flips and homological flops}
\author{Wai-Kit Yeung}
\address{Kavli IPMU, The University of Tokyo}
\email{wai-kit.yeung@ipmu.jp}
\begin{document}

\begin{abstract}
We give a survey for the results in \cite{Yeu20a, Yeu20b, Yeu20c}, which attempts to relate the derived categories under general classes of flips and flops. We indicate how the approach fails because of what appears to be a formal problem. We give some ideas, and record some failed attempts, to fix this problem. We also present some new examples.
\end{abstract}

\maketitle

%\vspace{-0.2cm}

\tableofcontents

\section{Log flips}

In this article, we give a survey for an approach developed in \cite{Yeu20a, Yeu20b, Yeu20c} to the problem of relating derived categories under flips and flops. The notions of flips and flops are sometimes a source of confusion because in the literature (especially those on derived categories) the terminology is sometimes used to refer to the closely related notions of $K$-dominant and $K$-equivalent birational maps. 
For the sake of clarity, we will recall the notions of flips and flops. In this section, we will start with the notion of a log flip (which generalizes both flips and flops), and we describe an algebraic setup associated to a log flip.
In the next section, we will define flips and flops, and elaborate more on this algebraic setup.
 
We do not assume knowledge of birational geometry from the reader. Most notions from birational geometry that we use will be reminded. We will also add several ``checkpoints'' in this article that summarize some essential information of the situation, so that readers less familiar with some parts of birational geometry or homological algebra will be able to read the rest of the article even if (s)he does not follow the arguments that lead to that description.

\bdf
A \emph{log flip} consists of birational projective morphisms between normal varieties
\begin{equation*}
	\begin{tikzcd} [row sep = 5]
		X^- \ar[rd, "\pi^-"'] & & X^+ \ar[ld, "\pi^+"] \\
		& Y
	\end{tikzcd}
\end{equation*}
such that both $\pi^-$ and $\pi^+$ are \emph{small}, meaning that $\Ex(\pi^{\pm}) \subset X^{\pm}$ have codimension $\geq 2$, 
together with Weil divisors $D^+, D^-, D^0$ on $X^+$, $X^-$, $Y$ respectively, strict transforms of each other, and satisfy
\begin{enumerate}
	\item $D^+$ is $\bQ$-Cartier and $\pi^+$-ample.
	\item $-D^-$ is $\bQ$-Cartier and $\pi^-$-ample.
\end{enumerate}
\edf

Recall that a Weil divisor $D \in \WDiv(X)$ is said to be $\bQ$-Cartier if $mD$ is Cartier for some $m \in \bZ\setminus \{0\}$. This allows one to discuss (relative) ampleness of $D$. By replacing them with a suitable multiple, one may take $D^{\pm}$ to be Cartier%
\footnote{Sometimes one would like to think of $D^{\pm}$ as part of the structure of a log flip, in which case one will have to tackle the case when $D^{\pm}$ are not Cartier.}. However, notice that $D^0$ is never Cartier (unless when $\pi^{\pm}$ are the identity maps).

Since $\pi^{\pm}$ are small, we have
\begin{equation}  \label{pushforward_O_DX}
	\pi^{\pm}_* \cO_{X^{\pm}}(iD^{\pm}) = \cO_Y(iD^0) \qquad \text{for any } i \in \bZ
\end{equation}

Let $\cA$ be the sheaf of $\bZ$-graded rings over $Y$ defined by $\cA = \bigoplus_{i \in \bZ} \cO_Y(iD^0)$, then by \eqref{pushforward_O_DX} and by relative ampleness of $D^+$ and $-D^-$, we have
\begin{equation}  \label{X_Proj_1}
	X^+ = \uProj^+_Y(\cA) := \uProj_Y(\cA_{\geq 0})
	\qquad \text{and} \qquad 
	X^- = \uProj^-_Y(\cA) := \uProj_Y(\cA_{\leq 0})
\end{equation}
Moreover, the reflexive sheaves $\cO_{X^{\pm}}(iD^{\pm})$ also admit a description in terms of $\cA$. Namely, denote by $\GrcA$ the category of quasi-coherent sheaves of graded $\cA$-modules. Then any $M \in \GrcA$ is also an object in $\Gr(\cA_{\geq 0})$, and hence determines an associated sheaf $\widetilde{M} \in \QCoh(\uProj^+_Y(\cA))$. The same is true for the negative direction. Then we have
\begin{equation*}
	\cO_{X^{\pm}}(iD^{\pm}) \cong \widetilde{\cA(i)} \qquad \text{under the identifications } \eqref{X_Proj_1}.
\end{equation*}

Everything we discuss here is local over $Y$, so we may assume that $Y$ is affine, in which case we summarize the essential information in the following
\begin{chkpt}  \label{chkpt_1}
Every log flip is locally described by a Noetherian $\bZ$-graded ring $A$. Namely, we have $R = A_0$, $Y = \Spec \, R$, $X^{\pm} = \Proj^{\pm}(A)$. Let $\cO_{X^{\pm}}(i) := \widetilde{A(i)}$, then we have $\pi^{\pm}_* (\cO_{X^{\pm}}(i)) = A_i$. 
\end{chkpt}

Here we mentioned that $A$ is Noetherian. It is relevant to our situation because of the following result (see, e.g., \cite[Theorem 1.5.5]{BH93}):
\bpp  \label{Noeth_gr_ring}
Let $A$ be a $\bZ$-graded ring. The the followings are equivalent:
\begin{enumerate}
	\item $A$ is a Noetherian ring; 
	\item every graded ideal of $A$ is finitely generated;
	\item $A_0$ is Noetherian, and both $A_{\geq 0}$ and $A_{\leq 0}$ are finitely generated over $A_0$;
	\item $A_0$ is Noetherian, and $A$ is finitely generated over $A_0$.
\end{enumerate}
\epp

We wish to study the derived category of $X^{\pm}$ in terms of the commutative algebra of $A$. The relation is obtained by Serre's equivalence. We have already seen above that there is a functor $(-)^{\sim} : \GrA \ra \QCoh(X^+)$. Let $\Torp \subset \GrA$ be the subcategory consisting of graded modules $M$ such that for every $x \in M$ there exists $s \in \bN$ such that $x \cdot (A_{>0})^s = 0$. Then the functor $(-)^{\sim}$ passes to the Serre quotient
\begin{equation}  \label{assoc_sheaf_quot_functor}
	(-)^{\sim} \, : \, \GrA / \Torp \raq \QCoh(X^+)
\end{equation}

The usual statement of Serre equivalence asserts that the functor \eqref{assoc_sheaf_quot_functor} is an equivalence if $A$ is non-negatively graded and generated in degree $0$ and $1$. If one sharpen each step of the usual proof, one have the following stronger version:

\bthm[\cite{Yeu20c}, Theorem 3.15]
Let $\cO_{X^+}(i) := \widetilde{A(i)}$ as above. Assume that the natural map
\begin{equation}  \label{pos_Cartier_map}
	\cO_{X^+}(i) \otimes_{\cO_{X^+}} \cO_{X^+}(j) \raq \cO_{X^+}(i+j)
\end{equation}
is an isomorphism for any $i,j \in \bZ$. Then the functor \eqref{assoc_sheaf_quot_functor} is an equivalence.
\ethm

In the case of log flips, the assumption of this Theorem is satisfied if and only if $D^+$ is Cartier. Motivated by this example, we will say that a Noetherian $\bZ$-graded ring $A$ is \emph{positively Cartier} if \eqref{pos_Cartier_map} is an isomorphism for any $i,j \in \bZ$.

At the level of derived categories, Serre's equivalence can be described in terms of a semi-orthogonal decomposition. More precisely, there is a local cohomology SOD (so called because of \eqref{local_cohom_seq}):
\begin{equation}  \label{local_cohom_SOD}
	\cD(\GrA) \, = \, \langle \, \cD_{\Iptr}(\GrA) \, , \, \cD_{\Torp}(\GrA)  \, \rangle
\end{equation}
Here, $\cD(\GrA)$ is the unbounded derived category of (not necessarily finitely generated) graded modules. 
$\cD_{\Torp}(\GrA) \subset \cD(\GrA)$ is the full triangulated subcategory consisting of objects $M \in \cD(\GrA)$ such that each cohomology module $H^i(M)$ is in the Serre subcategory $\Torp$. The inclusion $\iota : \cD_{\Torp}(\GrA) \rinto \cD(\GrA)$ has a right adjoint denoted as $\RGIp$, whose kernel is by definition $\cD_{\Iptr}(\GrA)$.

As the notation suggests, $\RGIp(M)$ is the local cohomology complex with respect to the graded ideal $I^+ := A_{>0} \cdot A$. If we denote the decomposition sequence of the SOD \eqref{local_cohom_SOD} by
\begin{equation}  \label{local_cohom_seq}
	\ldots \raq \RGIp(M) \xraq{\epsilon} M \xraq{\eta} \CIp(M) \raq \ldots 
\end{equation}
then $\CIp(M)$ can be described as a certain \v{C}ech complex, while $\RGIp(M)$ can be described as a certain extended \v{C}ech complex. Explicitly, if $I^+$ is generated by homogeneous elements $f_1,\ldots,f_r \in A_{>0}$, then they are given by
\begin{equation}  \label{local_cohom_Cech}
	\begin{split}
		\CIp(M) \, &:= \, \bigl[ \, \prod_{1 \leq i_0 \leq r} M_{f_{i_0}}  \raq 
		\prod_{1 \leq i_0 < i_1 \leq r} M_{f_{i_0}f_{i_1}}  \raq \ldots \raq
		M_{f_1\ldots f_r} \, \bigr] \\
		\RGIp(M) \, &:= \, \bigl[ \, M \raq \prod_{1 \leq i_0 \leq r} M_{f_{i_0}}  \raq
		\prod_{1 \leq i_0 < i_1 \leq r} M_{f_{i_0}f_{i_1}}  \raq \ldots \raq
		M_{f_1\ldots f_r} \, \bigr]
	\end{split}
\end{equation}
where the differentials are the alternating sums of restriction maps.

Here and below, if $X$ is a scheme or stack, then $\cD(X)$ will denote the unbounded derived category of quasi-coherent sheaves $\cD(X) := \cD(\QCoh(X))$.

If $A$ is positively Cartier, then it is a formal consequence of Serre's equivalence that the functor $(-)^{\sim} : \cD(\GrA) \ra \cD(X^+)$ restricts to an equivalence $\cD_{\Iptr}(\GrA) \ra \cD(X^+)$.
This is the basic setting we wish to exploit to pass from homological algebra on $\GrA$ to the derived category of $X^+$.
We summarize the situation in the following

\begin{chkpt}  \label{chkpt_2}
	For any Noetherian $\bZ$-graded ring $A$, there is a local cohomology SOD described by \eqref{local_cohom_SOD}, \eqref{local_cohom_seq}, \eqref{local_cohom_Cech}. If $A$ is positively Cartier then the semi-orthogonal component $\cD_{\Iptr}(\GrA)$ is equivalent to $\cD(X^+):= \cD(\QCoh(X^+))$.
	The same is true for the negative side.
\end{chkpt}

There is also a closely related stacky picture. While it will not be used in any technical way in this work, it will serve as a useful viewpoint.

Let $\frX$ be the stack $\frX = [\Spec (A) / \bG_m]$, so that $\GrA \simeq \QCoh(\frX)$. Let $\frX^{\pm} \subset \frX$ be the open substacks given by the complements of the closed subsets $V(I^{\pm}) \subset \Spec \, A$. Then $\frX^{\pm}$ are Deligne-Mumford, whose coarse moduli space are $X^{\pm}$. Serre's equivalence always holds with $X^{\pm}$ replaced by $\frX^{\pm}$. Moreover, $A$ is positively Cartier if and only if the map $\frX^{+} \ra X^{+}$ is an isomorphism.
Under this stacky viewpoint, the varieties $X^{\pm}$ (or rather their stacky versions $\frX^{\pm}$) are related by a variation of GIT quotients. 

\begin{chkpt}  \label{chkpt_3}
	The terms in the local cohomology sequence \eqref{local_cohom_seq} (at least for $M = A$) admit descriptions from two viewpoints
	\begin{center}
		\begin{tabular}{  c | c | c  }
		& 	$\pi^+ : X^+ \ra Y$  &  $j : \frX^+ \subset \frX $  \\ \hline
	$A$	& 	$A_i = \pi^+_* (\cO_{X^{+}}(i))$  &  $A = \cO_{\frX}$  \\ \hline
	$\CIp(A)$  & $\CIp(A)_i = {\bm R}\pi^+_* (\cO_{X^{+}}(i))$ & $\CIp(A) = {\bm R}j_* j^*(\cO_{\frX})$ \\	\hline
	$\RGIp(A)$  & $\RGIp(A)_i[1] = {\bm R}^{>0}\pi^+_* (\cO_{X^{+}}(i))$ & 
	\makecell{local cohomology along\\the unstable locus $V(I^+)$} 
		\end{tabular}
	\end{center}
%where we have written ``$=$'' for canonical isomorphisms under the identification $\cD(R) \simeq \cD(Y)$ and $\cD(\GrA) \simeq \cD(\frX)$ respectively.
%
Here and below, a subscript $(-)_i$ will always mean its weight grading $i$ part. Thus, given a complex $M$ of graded $A$-modules, then $M_i$ is a complex of $R$-module, where $R = A_0$. This gives a functor $(-)_i : \cD(\GrA) \ra \cD(R)$. This explains the meaning of the left column of this table. For example, the middle term means that $\CIp(A)_i \in \cD(R)$ and ${\bm R}\pi^+_* (\cO_{X^{+}}(i)) \in \cD(Y)$ are canonically identified under the equivalence $\cD(Y) \simeq \cD(R)$.
Indeed, this follows from the explicit \v{C}ech complex description \eqref{local_cohom_Cech}, whose weight $i$ part is exactly the same \v{C}ech complex%
\footnote{Since $X^+$ is defined as $\Proj(A_{\geq 0})$, the \v{C}ech complex that computes ${\bm R}\pi^+_* (\cO_{X^{+}}(i))$ is $\CIp(A_{\geq 0})_i$. However, for any homogeneous element $f$ with $\deg(f)>0$, we have $(A_{\geq 0})_f = A_f$, so that these two \v{C}ech complexes are the same. This is an instance of the fact that, although $X^+$ is defined as $\Proj(A_{\geq 0})$, it is really about $A$. Namely, it is covered by the affine open subschemes $\Spec (A_{(f)}) = \Spec ((A_{\geq 0})_{(f)})$ for $\deg(f) > 0$.} 
that computes ${\bm R}\pi^+_* (\cO_{X^{+}}(i))$. 
Notice that we have already seen in \eqref{chkpt_1} the first item of the left column. The third term is then a formal consequence of the first two
(here, the notation ${\bm R}^{>0}\pi^+_*$ is a shorthand for the good truncation $\tau^{>0} {\bm R}\pi^+_*$).

The equalities of the column on the right means that they are canonically identified under the equivalence $\cD(\GrA) \simeq \cD(\frX)$. Indeed, under this equivalence, the description \eqref{local_cohom_Cech} is standard.

\end{chkpt}

Generally speaking, when we discuss flips and flops in the next section, we will use the column on the left to extract homological information on $A$ from the situation of flips/flops, and then use the column on the right (as a viewpoint) to relate the derived categories of $X^-$ and $X^+$.

\section{Flips and flops}  \label{flip_flop_sec}

We consider flips and flops in this section. Before getting to that, we recall some results in algebraic geometry. The first is the following result, which may be found in \cite[Proposition 5.75]{KM98}:

\bthm  \label{OKX_omega}
Let $X$ be a projective normal variety%
\footnote{In this article, we will assume that $k$ is a perfect field. Hence a variety is smooth if and only if it is regular.}. Then there is an isomorphism $\cO(K_X) \cong \omega_X$, where $K_X$ is the canonical divisor, and $\omega_X$ is the dualizing sheaf.
\ethm
	
\brm  \label{canonical_dualizing_rmk}
This statement is still true if $X$ is not projective. Namely, for any variety $X$, take the canonical dualizing complex%
\footnote{Recall that for a morphism $f : X \ra Y$ between separated schemes of finite type over $k$, the functor $f^!$ is defined to be $f^! = j^* \circ g^!$, where $X \xra{j} \overline{X} \xra{g} Y$ is a factorization of $f$ such that $g$ is proper and $j$ is a Zariski open inclusion, and the functor $g^!$ is the right adjoint to ${\bm R}g_*$. All of these functors are taken between $\Dqcoh(-) \simeq \cD(\QCoh(-))$.} $\omega_X^{\bullet} := p^!(k)$, where $p : X \ra \Spec \, k$ is the projection to a point, and define $\omega_X := \cH^{-n}(\omega_X^{\bullet})$ where $n = \dim(X)$. If $X$ is proper (over $k$), then $\omega_X$ can be characterized as the dualizing sheaf. But even if $X$ is not proper,  $\omega_X \in \Coh(X)$ is still canonically defined. Moreover, if $X$ is normal then we have $\cO(K_X) \cong \omega_X$, where $K_X$ is the canonical divisor. For example, if $X$ is quasi-projective then this follows from Theorem \ref{OKX_omega} by taking a compactification. This implies that $\omega_X$ is reflexive whenever $X$ is a normal variety (since we can argue on each affine open subscheme, which is quasi-projective). Since $\omega_X$ and $\cO(K_X)$ are both reflexive and coincide on the smooth locus, they must be isomorphic.
\erm

We will also need a version%
\footnote{See \cite[Theorem 1-2-3]{KKM87} for the case when $X$ is smooth. This implies our present case by taking a resolution of singularities $f : \widetilde{X} \ra X$. Namely, since $X$ has rational singularities, we have ${\bm R}f_* (\omega_{\widetilde{X}}) \cong \omega_X$, so that we have ${\bm R}f_* (\omega_{\widetilde{X}} \otimes f^* \scL) \cong \omega_X \otimes \scL$ by the projection formula. Then apply ${\bm R}(\pi \circ f)_* \cong {\bm R}\pi _* \circ {\bm R}f_*$ to $\omega_{\widetilde{X}} \otimes f^* \scL$.} of the Kawamata-Viehweg vanishing theorem:

\bthm  \label{KV_vanishing}
Assume ${\rm char}(k)=0$.
Let $\pi : X \ra Y$ be a projective birational morphism between varieties, where $X$ has at most rational singularities. Let $\scL \in \Pic(X)$ be $\pi$-nef, then we have ${\bm R}^i\pi_*(\omega_X \otimes \scL) = 0$ for all $i > 0$.
\ethm

Now we recall the notion of flips and flops, which are two mutually exclusive subclasses of log flips.
We will say that a log flip is trivial if both $\pi^{\pm}$ are the identity maps. Otherwise, we say that the log flip is non-trivial.

\bdf
A \emph{flip} is a non-trivial log flip $(X^{\pm},Y,\pi^{\pm},D^{\pm})$ in which $D^{\pm} = K_{X^{\pm}}$.

A \emph{flop} is a non-trivial log flip in which $K_{X^{\pm}}$ is numerically $\pi^{\pm}$-trivial.

A \emph{Gorenstein flop} is a non-trivial log flip in which $K_Y$ is Cartier. By smallness, we then have $(\pi^{\pm})^* K_Y = K_{X^{\pm}}$, and therefore it is a flop.
\edf

We will consider flips and Gorenstein flops. In fact, we will impose one more condition:
\begin{equation}  \label{rat_Gor_assump}
	\parbox{40em}{Assume that $X^{\pm}$ have at most rational Gorenstein singularities.}
\end{equation}
In the case of flips or Gorenstein flops, by applying Theorem \ref{KV_vanishing} to $\pi^- : X^- \ra Y$, the assumption \eqref{rat_Gor_assump} implies that $Y$ has at most rational singularities, and hence is Cohen-Macaulay. In this case, since $X^{\pm}$ and $Y$ are Cohen-Macaualay, their canonical dualizing complexes (see Remark \ref{canonical_dualizing_rmk}) is concentrated in degree $-n$, so that $\omega_X^{\bullet} \cong \omega_X[n]$ for $X$ being $X^{\pm}$ or $Y$.

Thus, in the case of flips satisfying \eqref{rat_Gor_assump}, we have
\begin{equation*}
	(\pi^{\pm})^!(\omega_Y) \, = \, \omega_{X^{\pm}} \, \stackrel{\eqref{OKX_omega}}{\cong} \, \cO(K_{X^{\pm}}) \, = \, \cO_{X^{\pm}}(1)
\end{equation*}
(the reader is cautioned that $\cO_{X^{-}}(1)$ is $\pi^-$-anti-ample in our convention).

In the case of Gorenstein flops satisfying \eqref{rat_Gor_assump}, we have
\begin{equation*}
	(\pi^{\pm})^!(\cO_Y) \, \cong \,  \cO_{X^{\pm}}
\end{equation*}
Indeed, if $f : X \ra Y$ is a proper birational morphism between Gorenstein normal varieties such that $f^* K_Y = K_X$, then we have $f^*(\omega_Y) \cong \omega_X$ by \eqref{OKX_omega}. Also, we again have $f^!(\omega_Y) \cong \omega_X$ by virtue of $X$ and $Y$ being Cohen-Macaulay. Since $f^!(\cF) \cong f^*(\cF) \otimes f^!(\cO_Y)$ for perfect complexes $\cF \in \Dperf(Y)$, we concludes that $f^!(\cO_Y) \cong \cO_X$.

Thus, in both cases, we have the following

\begin{chkpt}  \label{chkpt_4}
	For flips and Gorenstein flops satisfying \eqref{rat_Gor_assump}, there is a dualizing complex $\omega_Y^{\bullet}$ on $Y$ such that
	\begin{equation*}
		(\pi^{\pm})^!(\omega_Y^{\bullet}) \, \cong \, \cO_{X^{\pm}}(a)
	\end{equation*}
where $a=0$ for  Gorenstein flops and $a=1$ for flips.
\end{chkpt}

Our goal now is to translate this condition into $\cD(\GrA)$, in the spirit of the last section. First, we introduce the following degreewise dualizing functor%
\footnote{Of course, the weight-by-weight description \eqref{D_Y_def} of $\bD_Y(M)$ is not enough to determine $\bD_Y(M)$ as an object in $\cD(\GrA)$. To actually define $\bD_Y$, one could start with $\Homcom_R(-,-):\GrA^{\op} \times \Mod(R) \ra \GrA$ and take its derived functor.} (recall that $R = A_0$ and $Y = \Spec \, R$).
\begin{equation}  \label{D_Y_def}
	\bD_Y : \cD(\GrA)^{\op} \ra \cD(\GrA) \, , \quad  \qquad \bD_Y(M)_i = \RHom_R(M_{-i},\omega_Y^{\bullet})
\end{equation}

We postulate the following as an imitation of \eqref{chkpt_4}:
\bdf
A \emph{homological flip/flop} consists of a pair of isomorphisms in $\cD(\GrA)$:
\begin{equation}  \label{Phi_maps}
	\begin{split}
		\Phi^+ &: \CIp(A)(a) \xra{\cong} \bD_Y( \CIp(A) ) \\
		\Phi^- &: \CIm(A)(a) \xra{\cong} \bD_Y( \CIm(A) ) 
	\end{split}
\end{equation}
such that they are compatible, in the sense that the following diagram commutes:
	\begin{equation}  \label{Phi_compatible}
	\begin{tikzcd} [row sep = 0]
		& \CIp(A)(a) \ar[r, "\Phi^+", "\cong"'] & \bD_Y( \CIp(A) ) \ar[rd, "\bD_Y(\eta^+)"] & \\
		A(a) \ar[ru, "\eta^+"] \ar[rd, "\eta^-"']  & & & \bD_Y(A) \\
		& \CIm(A)(a) \ar[r, "\Phi^-", "\cong"'] & \bD_Y( \CIm(A) ) \ar[ru, "\bD_Y(\eta^-)"'] &
	\end{tikzcd}
\end{equation}
In the case $a = 0$, we call it a homological flop. In the case $a = 1$, we call it a homological flip%
\footnote{One may call the general notion a homological flap, so that it becomes a fl0p when $a=0$, and a fl1p when $a = 1$.}.
\edf

The existence of the isomorphisms $\Phi^{\pm}$ should be quite believable if we look at both sides of \eqref{Phi_maps} weight by weight. For example, by Grothendieck duality, we have 
\begin{equation}  \label{Phi_isom_weight_i}
	\RcHom_{\cO_{Y}}( {\bm R}\pi^{+}_* (\cO_{X^{+}}(-i))  , \omega_Y^{\bullet} )\, \cong \, {\bm R}\pi^{+}_* \RcHom_{\cO_{X^{+}}}( \cO_{X^{+}}(-i) , \pi^!( \omega_Y^{\bullet} ) )  \stackrel{\eqref{chkpt_4}}{\cong}  {\bm R}\pi^{+}_* (\cO_{X^{+}}(i+a))
\end{equation}
In view of \eqref{chkpt_3}, the left and right hand side are $\bD_Y(\CIp(A))_i$ and $\CIp(A)_{i+a}$ respectively.

The meaning of the commutative diagram \eqref{Phi_compatible} is harder to describe. Let us just say that it boils down to the fact that the isomorphism $\Phi_X : \cO(K_X) \xra{\cong} \omega_X$ in Theorem \ref{OKX_omega} can be appropriately chosen so that they satisfy a certain compatibility condition across the maps $\pi^{\pm}$. More precisely, we may require the following diagram to commute:
\begin{equation}  \label{Phi_maps_commute}
	\begin{tikzcd}
		\pi^-_* \cO(K_{X^-})  \ar[d, "\pi^-_*(\Phi_{X^-})"'] \ar[rr, equal]
		& & \cO(K_{Y})  \ar[d, "\Phi_{Y}"]
		& & \pi^+_* \cO(K_{X^+}) \ar[d, "\pi^+_*(\Phi_{X^+})"] \ar[ll, equal]  \\
		\pi^-_* \omega_{X^-} \ar[rr, "\Tr_{X^-/Y}"]
		& &  \omega_{Y}
		& & \pi^+_* \omega_{X^+}  \ar[ll, " \Tr_{X^+/Y}"']
	\end{tikzcd}
\end{equation}
where $\Phi_{X^{\pm}}$ and $\Phi_Y$ are the isomorphisms in Theorem \ref{OKX_omega}, and $\Tr_{X^{\pm}/Y}$ are induced from certain adjunction morphisms on dualizing complexes. Here, we have written equalities in the first row because they are equal as subsheaves of the sheaf of rational functions on $Y$.

In fact, to pass from the weight-by-weight description \eqref{Phi_isom_weight_i} to an actual isomorphism \eqref{Phi_maps} in $\cD(\GrA)$, we had to develop from scratch the effect of Grothendieck duality in $\cD(\GrA)$ (this is performed in \cite{Yeu20c}). Then, we had to show that, under this formalism, \eqref{Phi_maps_commute} indeed translates into \eqref{Phi_compatible}. For this, we had to check many commutative diagrams. All these checkings are lengthy and tedious (see \cite{Yeu20a} for details), but doesn't seem to offer any additional insight, so we will skip it in this article. We summarize this discussion by the following
\bthm  \label{satisfies_HPP}
For flips and Gorenstein flops satisfying \eqref{rat_Gor_assump}, the corresponding%
\footnote{For flips, there is no ambiguity in choosing $D^{\pm}$, so that $A$ is canonically defined. For flops, we will assume that $D^{\pm}$ are chosen to be Cartier, possibly by replacing them by their multiple. Since we work under the assumption \eqref{rat_Gor_assump}, we see that $D^{\pm}$ are Cartier in both cases. In fact, this was already implicitly used in the second isomorphism in \eqref{Phi_isom_weight_i}, which holds because $\cO_{X^{+}}(-i)$ is locally free. Notice that it also holds if $\cO_{X^{+}}(-i) = \cO_{X^{+}}(-iD^+)$ is maximal Cohen-Macaulay, which is satisfied whenever $X^{\pm}$ have at most log terminal singularities (see \cite[Corollary 5.25]{KM98}). Accordingly, the assumption \eqref{rat_Gor_assump} can be weakened.} graded ring $A$ admits the structure of a homological flip/flop.
\ethm

We will formulate one more property of $A$ that imitates the Kawamata-Viehweg vanishing theorem. 
We will write $\cD_{<w}(\GrA)$ as the full subcategory of $\cD(\GrA)$ consisting of $M \in \cD(\GrA)$ such that $M_i = 0$ for all $i \geq w$. Again, we emphasize here that $M_i$ refers to the weight degree, not the homological degree. In particular, $\cD_{<w}(\GrA)$ is a triangulated subcategory.
The subcategory $\cD_{>w}(\GrA)$ is also defined in a similar way.

\bdf
A $\bZ$-graded ring $A$ is said to satisfy \emph{canonical vanishing} at $a \in \bZ$ if we have $\RGIp(A) \in \cD_{<a}(\GrA)$ and $\RGIm(A) \in \cD_{>a}(\GrA)$.
\edf

In view of the table of \eqref{chkpt_3}, we see that Theorem \ref{KV_vanishing} immediately implies the following

\bpp  \label{satisfies_can_van}
For flips and Gorenstein flops satisfying \eqref{rat_Gor_assump}, the corresponding graded ring $A$ satisfies canonical vanishing at $a$ (where $a = 0$ for flops and $a = 1$ for flips).
\epp

Next, we explore some consequences of Theorem \ref{satisfies_HPP} and Proposition \ref{satisfies_can_van}. The following Theorem, as well as its proof, may be regarded as the main result of this work:
\bthm  \label{main_thm_1}
Suppose $(A,\Phi^+,\Phi^-)$ is a homological flip/flop that satisfies canonical vanishing at $a$ (where $a = 0$ for flops and $a = 1$ for flips), then there is an isomorphism in $\cD(\GrA)$:
\begin{equation*}
\Psi \, : \, \RGIp(A)(a)[1] \xraq{\cong} \bD_Y(\RGIm(A))
\end{equation*}
Moreover, $A$ is Gorenstein.
\ethm

\bpf
Consider the commutative diagram \eqref{Phi_compatible}. Since the maps $\Phi^{\pm}$ are isomorphisms, we may regard it as a commutative square. As such, the $3 \times 3$-lemma (see, e.g., \cite[Proposition 1.1.11]{BBD82} or \cite[Lemma 2.6]{May01}) asserts that it can be extended to a $3 \times 3$ square. More precisely, the object $Z$ as well as the maps in the dotted lines of the followng diagram exists, making each row and column part of a distinguished triangle:
\begin{equation}  \label{HFF_big_diag_1}
	\begin{tikzcd}
		A(a) \ar[r, "\eta^+"] \ar[d, "\eta^-"] & \CIp(A)(a) \ar[d, "\bD_Y(\eta^+) \circ \Phi^+"] \ar[r, "\delta^+"] 
		& \RGIp(A)(a)[1] \ar[d, dashed, "\Psi"] \\
		\CIm(A)(a) \ar[r, "\bD_Y(\eta^-) \circ \Phi^-"]  \ar[d, "\delta^-"]
		& \bD_Y(A) \ar[r, "\bD_Y(\epsilon^-)"] \ar[d, "\bD_Y(\epsilon^+)"] 
		& \bD_Y(\RGIm(A)) \ar[d, dashed]\\
		\RGIm(A)(a)[1] \ar[r, dashed, "\Psi'"] &  \bD_Y(\RGIp(A)) \ar[r, dashed] & Z
	\end{tikzcd}
\end{equation}

By the assumption of canonical vanishing, we have $\RGIp(\cA)(a)[1] \in \cD_{<0}(\GrA)$ and $\bD_Y(\RGIm(\cA)) \in \cD_{<-a}(\GrA)$, so that, as the cone of $\Psi$, we have $Z \in \cD_{<0}(\GrA)$. Likewise, we have $\RGIm(\cA)(a)[1] \in \cD_{>0}(\GrA)$ and $\bD_Y(\RGIp(\cA)) \in \cD_{>-a}(\GrA)$, so that as the cone of $\Psi'$, we have $Z \in \cD_{>-a}(\GrA)$. Since $a = 0$ or $a = 1$, we must have $Z = 0$. {\it i.e.,} $\Psi$ is an isomorphism.

To prove that $A$ is Gorenstein, we will show that it has finite injective dimension%
\footnote{Given a graded module $M \in \GrA$, then it has finite injective dimension in the abelian category $\GrA$ if and only if it has finite injective dimension in the abelian category $\Mod(A)$.
In fact, we have ${\rm inj \, dim}_{\GrA}(M) \, \leq \, {\rm inj \, dim}_{\Mod(A)}(M) \, \leq \, {\rm inj \, dim}_{\GrA}(M) + 1$ (see, {\it e.g.} \cite[Proposition 3.6.6]{BH93}). Hence there is no ambiguity in our discussion.}.
Take the local cohomology sequence
\begin{equation*}  
	\ldots \raq \RGIp(A) \raq A \raq \CIp(A) \raq \ldots 
\end{equation*}
Then the terms can be rewritten as $\RGIp(A) \cong \bD_Y( \RGIm(A) )(-a)[-1]$ and $\CIp(A) \cong \bD_Y( \CIp(A) )(-a)$. Notice that both $\RGIm(A)$ and $\CIp(A)$ have finite Tor dimension because of the explicit presentation \eqref{local_cohom_Cech}. Therefore their $\bD_Y$-dual have finite injective dimension, hence so does $A$.
\epf

For general classes of flips and flops, it seems quite difficult to describe this duality between local cohomology groups explicitly. There is however one simple example where such a description is possible:

\beg  \label{polynomial_RGam_example}
Let $A = k[x_1,\ldots,x_p,y_1,\ldots,y_q]$, where $\deg(x_i)=1$ and $\deg(y_i)=-1$. This corresponds to the standard flip/flop, and we have $a = q-p$. By the explicit presentation \eqref{local_cohom_Cech}, it is easy to see that
\begin{equation*}
	\RGIp(A) \, = \, k[y_1,\ldots,y_q] \otimes \RGIp(k[x_1,\ldots,x_p])
\end{equation*}
Also, in view of the table in \eqref{chkpt_3} again, we see that $\RGIp(k[x_1,\ldots,x_p])_i$ is simply the higher cohomology of $\bP^{p-1}$ (shifted by $1$) of $\cO(i)$. Such a computation is standard (see, {\it e.g.,} \cite[Section III.5]{Har77}), and we have 
\begin{equation*}
	\RGIp(A) \, = \, k[y_1,\ldots,y_q] \otimes k[x_1,\ldots,x_p]^* (p)[-p]
\end{equation*}
where we denote $M^*$ to be the $k$-linear dual of $M$. Similarly, we have
\begin{equation*}
	\RGIm(A) \, = \, k[x_1,\ldots,x_p] \otimes k[y_1,\ldots,y_q]^* (-q)[-q]
\end{equation*}

We want to claim that they are $\bD_Y$-dual to each other. While the description of $\bD_Y$ on $R = A_0$ seems to be quite complicated in general, we will only need a special case. Namely, since $\RGIp(A)_i$ corresponds to higher pushforwards of $\cO_{X^{+}}(i)$, it must be supported on the image of the exceptional locus, which is a point. In this case, $\bD_Y$ is easy to describe:
\blm
Suppose that $M \in \Dbcoh(R)$ is (set-theoretically) supported on $\Spec \, k = \Spec \, R/J$, and let $\omega_{Y}^{\bullet} \in \Dbcoh(R)$ be a dualizing complex, normalized so that its cohomology sheaf is concentrated%
\footnote{$R$ is Cohen-Macaulay in this example, so that we may assume that $\cH^{\bullet}(\omega_{Y}^{\bullet})$ is concentrated in degree $0$. An analogous statement holds for more general $R$ (at least for domains finitely generated over $k$), for which we may assume that $\cH^{\bullet}(\omega_{Y}^{\bullet})$ is concentrated in degree $0$ on the smooth part.} 
in degree $0$. Then we have $\bD_Y(M) \cong M^*[-n]$, where $n = \dim R = p+q-1$.
\elm

%\bpf
%Write $M = i_*(N)$ for some $i : Z := \Spec (R/J^n) \ra \Spec \, R$, then we have $\bD_Y(i_*N) = i_* \bD_Z(N)$, where $\bD_Z$ is the dualizing functor defined by the dualizing complex $i^!(\omega_{Y}^{\bullet})$. Since $R/J^n$ is a local ring, it has only one dualizing complex up to homological shifts. Also, since $R/J^n$ is $0$-dimensional, $(-)^*$ is a dualizing functor, and hence must coincide with $\bD_Z$ up to a shift, which is determined from the normalization of $\omega_{Y}^{\bullet}$.
%\epf

From this, we see that
\begin{equation*}
		\bD_Y(\RGIm(A)) \, = \, k[x_1,\ldots,x_p]^* \otimes k[y_1,\ldots,y_q] (q)[q][-p-q+1)] \, = \, \RGIp(A)(a)[1]
\end{equation*}
\eeg

\brm  \label{poly_ring_RGam_rmk}
This example can be directly generalized to the case $A = k[x_1,\ldots,x_p,y_1,\ldots,y_q]$, where $\deg(x_i)>0$ and $\deg(y_i)<0$.
In this case, let $\eta^+ = \sum_{i=1}^p \deg(x_i)$ and $\eta^- = -\sum_{j=1}^p \deg(y_j)$, and take $a = \eta^- - \eta^+$.
Notice that $\RGIp(A)$ is still the same extended \v{C}ech complex as in Example \ref{polynomial_RGam_example}, but with generators in different degrees. As such, the computation in \cite[Section III.5]{Har77} carries through verbatim, and we have
\begin{equation*}
	\begin{split}
		\RGIp(A) \, &= \, k[y_1,\ldots,y_q] \otimes k[x_1,\ldots,x_p]^* (\eta^+)[-p] \\
		\RGIm(A) \, &= \, k[x_1,\ldots,x_p] \otimes k[y_1,\ldots,y_q]^* (-\eta^-)[-q]
	\end{split}
\end{equation*}
and hence we still have $\bD_Y(\RGIm(A)) \cong \RGIp(A)(a)[1]$.
\erm

We wish to apply the results in Theorem \ref{satisfies_HPP}, Proposition \ref{satisfies_can_van} and Theorem \ref{main_thm_1} to relate the derived categories of $X^-$ and $X^+$. We may summarize the situation as follows
\begin{equation*}
	\begin{tikzcd} [row sep = 5]
		& \cD(\GrA) \ar[ld,"\text{kill }\RGIm"']  \ar[rd, "\text{kill }\RGIp"]& \\
		\cD(X^-) & & \cD(X^+)
	\end{tikzcd}
\end{equation*}
where the derived categories $\cD(X^{\pm})$ are obtained by ``killing'' the endofunctors ${\bm R}\Gamma_{I^{\pm}}(M) = M \otimes_A^{{\bm L}} {\bm R}\Gamma_{I^{\pm}}(A)$.
The duality of local cohomology groups in Theorem \ref{main_thm_1} then suggests a way to relate these two derived categories. However, an implementation of this idea seems to be not so straightforward. We present one attempt in the next section. This approach is successful in some useful cases (see Section \ref{example_sec}), but fall short in general because of what seems to be a formal problem.

\section{Weight truncation}  \label{wtr_sec}

Before describing the actual construction, we would like to suggest a parallelism with a paper \cite{Orl09} of Orlov. Namely, we may summarize our desired picture as follows:
\begin{equation}  \label{expected_pic_flip_flop}
	\begin{split}
		&\text{Flip or Gorenstein flop (with good singularities)} \\
	\Longrightarrow	& \, A \text{ is Gorenstein, and } \RGIp(A)(a)[1] \cong \bD_Y(\RGIm(A)) \\
	\stackrel{(?)}{\Longrightarrow}	& \, 
	\begin{cases}
		\Dperf(X^-) \, \hookleftarrow \, \Dperf(X^+) & \text{ if } a > 0,\\
		\Dperf(X^-) \, \simeq \, \Dperf(X^+) & \text{ if } a = 0,\\
		\Dperf(X^-) \, \hookrightarrow \, \Dperf(X^+) & \text{ if } a < 0
	\end{cases}
	\end{split}
\end{equation}

On the other hand, for a Gorenstein projective variety $X$ of dimension $n$, Orlov studied the relation between $\Dbcoh(X)$ and the triangulated category of singularities of the projective cone of $X$. More precisely, assume that $\cO(1)$ is very ample, and let $A = \bigoplus_{i \geq 0} H^0(X,\cO(i))$. Assume that $H^j(X,\cO(i)) = 0$ for all $j \neq 0,n$ for all $i \in \bZ$. Then Orlov established in \cite{Orl09} the following implications
\begin{equation}  \label{Orlov_picture}
	\begin{split}
		& \omega_X \cong \cO(-a) \\
		\Longrightarrow	& \, A \text{ is Gorenstein, and } \RHomcom_A(k,A) \cong k(a)[-n] \\
		\stackrel{(*)}{\Longrightarrow}	& \, 
		\begin{cases}
			\cD_{{\rm sg}}({\rm gr}(A)) \, \hookrightarrow \, \Dbcoh(X) & \text{ if } a > 0,\\
			\cD_{{\rm sg}}({\rm gr}(A)) \, \simeq \, \Dbcoh(X) & \text{ if } a = 0,\\
			\cD_{{\rm sg}}({\rm gr}(A)) \, \hookleftarrow \, \Dbcoh(X) & \text{ if } a < 0
		\end{cases}
	\end{split}
\end{equation}

Our goal in this section is to imitate the construction of the step $(*)$ in \eqref{Orlov_picture} and try to establish the step $(?)$ in \eqref{expected_pic_flip_flop}. Our construction is also influenced and motivated by the papers \cite{HL15,BFK19}.

Fix an integer $w \in \bZ$ once and for all. Let $\cD_{[\geq w]}(\GrA)$ be the smallest cocomplete ({\it i.e.,} closed under arbitrary direct sums) triangulated subcategory of $\cD(\GrA)$ containing the objects $\{A(-i)\}_{i \geq w}$, and let $\cD_{<w}(\GrA)$ be defined as above. {\it i.e.,} it consists of $M \in \cD(\GrA)$ such that $M_i = 0$ for all $i \geq w$. By Neeman-Brown representability, we have an SOD
\begin{equation}  \label{weight_trunc_SOD}
	\cD(\GrA) \, = \, \langle \,  \cD_{<w}(\GrA) \, , \,  \cD_{[\geq w]}(\GrA)  \, \rangle
\end{equation}

The terminology $\cD_{[\geq w]}(\GrA)$ refers to those objects that are generated in weight $\geq w$, while $\cD_{<w}(\GrA)$ refers to those that are concentrated in weight $<w$. In the setting of \cite{Orl09}, $A$ is non-negatively graded, so that $\cD_{[\geq w]}(\GrA) = \cD_{\geq w}(\GrA)$. But in our setting of $\bZ$-graded rings, these two are very different. In particular, none contain the other.

The decomposition sequence associated to \eqref{weight_trunc_SOD} will be denoted as
\begin{equation}  \label{weight_trunc_seq}
	\ldots \raq \cL_{[\geq w]}M \raq M \raq \cL_{<w}M \raq \ldots
\end{equation}

We will put together the SODs \eqref{local_cohom_SOD} and \eqref{weight_trunc_SOD}. Notice that $\cD_{<w}(\GrA)$ consists of those $M \in \cD(\GrA)$ such that each $H^i(M)$ has weight concentrated in degree $<w$, and hence we have $\cD_{<w}(\GrA) \subset \cD_{\Torp}(\GrA)$. Thus, if we compare \eqref{local_cohom_SOD} and \eqref{weight_trunc_SOD}, we would expect that $\cD_{[\geq w]}(\GrA)$ would be ``bigger'' than $\cD_{\Iptr}(\GrA)$. In fact, one can show that $\cD_{\Iptr}(\GrA)$ embeds (via a non-identity functor) as a semi-orthogonal summand of $\cD_{[\geq w]}(\GrA)$. Combined with \eqref{weight_trunc_SOD}, it then gives us a $3$-term semi-orthogonal decomposition of $\cD(\GrA)$. This is summarized in the following result, whose proof is simple once it is formulated precisely as below:

\bthm
The restriction of $\cL_{[\geq w]}$ to $\cD_{\Iptr}(\GrA)$ gives a fully faithful functor $\cL_{[\geq w]} : \cD_{\Iptr}(\GrA) \ra  \cD_{[\geq w]}(\GrA)$, with a left adjoint given by $\CIp$. Hence there is a semi-orthogonal decomposition 
\begin{equation*}
\cD_{[\geq w]}(\GrA) \, = \,	\langle \, \cL_{[\geq w]} \cD_{\Iptr}(\GrA)  \, , \, \cD_{[\geq w], \Torp}(\GrA)   \, \rangle
\end{equation*}

Combined with \eqref{weight_trunc_SOD}, there is therefore a $3$-term semi-orthogonal decomposition
\begin{equation}  \label{three_term_SOD}
	\cD(\GrA) \, = \, \langle \, \UOLunderbrace{ \cD_{<w}(\GrA) \, , \,} 
	[\cL_{[\geq w]} \cD_{\Iptr}(\GrA)]_{=: \, \cD_{(<w)}(\GrA)} \UOLoverbrace{  \, , \, \cD_{[\geq w], \Torp}(\GrA) }^{\cD_{[\geq w]}(\GrA)}  \, \rangle
\end{equation}
where the subcategory $\cD_{(<w)}(\GrA) \subset \cD(\GrA)$ can be characterized more precisely as
\begin{equation*}
	\cD_{(<w)}(\GrA) \, = \, \{ \, M \in \cD(\GrA) \, | \, \RGIp(M) \in \cD_{<w}(\GrA) \, \}
\end{equation*}

The middle component, called the window subcategory can therefore be characterized as
\begin{equation*}
	\cL_{[\geq w]} \cD_{\Iptr}(\GrA) \, = \, \cD_{[\geq w], (<w)}(\GrA) \, = \, \cD_{[\geq w]}(\GrA) \cap \cD_{(<w)}(\GrA)
\end{equation*}
which is equivalent to $\cD_{\Iptr}(\GrA)$ (and hence to $\cD(\frX^{+})$) via the functors $\CIp$ and $\cL_{[\geq w]}$.
\ethm

Notice that, in this Theorem, we have embedded $\cD(\frX^{+})$ into $\cD(\GrA)$ in a non-standard way, {\it i.e.,} as the window subcategory $\cD_{[\geq w], (<w)}(\GrA)$, instead of the more straightforward $\cD_{\Iptr}(\GrA)$. Notice that the straightforward embedding to $\cD_{\Iptr}(\GrA)$ has the disadvantage that $\CIp(M)$ almost never have coherent cohomology even if $M$ does (it has finite cohomological dimension though!). In contrast, this non-standard ``window embedding'' have the following advantage:

\bthm
The $3$-term SOD \eqref{three_term_SOD} restrict to a $3$-term SOD on $\Dmcoh(\GrA)$:
\begin{equation*}
	\Dmcoh(\GrA) \, = \, \langle \, \cD^-_{{\rm coh},<w}(\GrA) \, , \,  \cD^-_{{\rm coh},[\geq w],(<w)}(\GrA)  \, , \,  \cD^-_{{\rm coh},[\geq w],\Torp}(\GrA)  \, \rangle
\end{equation*}
where the middle component is equivalent to $\Dmcoh(\frX^+)$.
\ethm

%The interested reader may find a proof in \cite{aaaaa}, although we advise him/her to wait for a revised version, where the arguments will be simplified and put in a broader scope.

A proof may be found in \cite{Yeu20b} (a revised version is under preparation where the arguments will be simplified and put in a broader scope), and is based on the tensor form \eqref{cL_tensor_form} below. We will skip the proof in this article.

Now we are ready to describe the functor that relate the derived categories under flips/flops. First, we introduce the functor
\begin{equation*}
	\bD_A \, : \, \cD(\GrA)^{\op} \raq \cD(\GrA) \, , \qquad \bD_A(M) := \RHomcom_A(M,A)
\end{equation*}
We postulate the following idea:
\begin{equation}  \label{DA_idea_1}
	\parbox{40em}{The functor $\bD_A$ ``should'' send the window $\cD_{[\geq w], (<w)}(\GrA)$ in the positive direction to the window $\cD_{[\leq -w], (>-w)}(\GrA)$ in the negative direction.}
\end{equation}

The idea does not work as it is stated. But let's press on with the idea, which may be broken down into two parts:
\begin{equation}  \label{DA_idea_2}
	\parbox{40em}{
		\begin{enumerate}
			\item The functor $\bD_A$ ``should'' send $\cD_{[\geq w]}(\GrA)$ to $\cD_{[\leq -w]}(\GrA)$
			\item The functor $\bD_A$ ``should'' send $\cD_{(<w)}(\GrA)$ to $\cD_{(>-w)}(\GrA)$
		\end{enumerate}
		}
\end{equation}

The first item seems quite reasonable, because $\bD_A$ sends the generators $\{A(-i)\}_{i \geq w}$ of $\cD_{[\geq w]}(\GrA)$ to the generators $\{A(-i)\}_{i \leq -w}$ of $\cD_{[\leq -w]}(\GrA)$. However, notice that $\bD_A$ sends infinite direct sums to infinite direct products, and $\cD_{[\leq -w]}(\GrA)$ seems to be not closed under infinite direct products in general, so that (1) seems to be not true in general. In fact, it seems to be not true even if one restrict to $\cD^-_{{\rm coh},[\geq w]}(\GrA)$. More precisely, one may encounter the following problem:
\begin{equation}  \label{bdd_below_complex_weight}
	\parbox{40em}{Suppose $M$ is a bounded below complex $M = [0 \ra M^n \ra M^{n+1} \ra \ldots]$ such that each $M^j$ is a finite direct sums of the free graded modules $A(-i)$ for $i \leq -w$, then it is not clear whether $M$ is in $\cD_{[\leq -w]}(\GrA)$.}
\end{equation}

In contrast, one can actually prove the item (2), when restricted to $\Dmcoh(\GrA)$:
\bpp
Suppose there is a dualizing complex $\omega_Y^{\bullet}$ on $R = A_0$ such that there is an isomorphism $\RGIp(A)(a)[1] \cong \bD_Y(\RGIm(A))$ in $\cD(\GrA)$. Then $\bD_A$ sends $\cD^-_{{\rm coh},(<w)}(\GrA)$ to $\cD^+_{{\rm coh},(>-w+a)}(\GrA)$.
\epp

\bpf
Suppose $M \in \Dmcoh(\GrA)$ satisfies $\RGIp(M) \in \cD_{<w}(\GrA)$, then we have
\begin{equation*}
	\begin{split}
		\RGIm(\bD_A(M)) \, &= \, \RHomcom_A(M,A) \otimes_A^{{\bm L}} \RGIm(A) \\
		\, &\stackrel{(*)}{=} \, \RHomcom_A(M,\RGIm(A))  \\
		 \, &\cong \,  \bD_Y(M \otimes_A^{{\bm L}} \RGIp(A))(-a)[-1]
	\end{split}
\end{equation*}
where $(*)$ holds because $\RGIm(A)$ has finite Tor dimension and $M \in \Dmcoh(\GrA)$.
\epf

By the results of the last section, we may assume that $A$ is Gorenstein, so that $\bD_A$ is in some sense well-behaved. Thus, if \eqref{DA_idea_2}(1) holds in some good context, we would have obtained a relation between the derived categories of $\frX^+$ and $\frX^-$ similar to the expected relation \eqref{expected_pic_flip_flop}. 
While this picture of using $\bD_A$ to relate the two windows was the first argument along this line that the author discovered (and is still very attractive to the author), it seems to be not the most efficient one. We now seek to reformulate it without using the duality $\bD_A$.

%\brm  \label{formal_skip_rmk}
%In the rest of this section, we will give some formal reformulations of the idea \eqref{DA_idea_1}. The author has spent/wasted a lot of time trying to find a good formulation so that problem like \eqref{bdd_below_complex_weight} will not affect the argument. It seems that in each of these reformulation, problems like \eqref{bdd_below_complex_weight} invariably show up in a disguised form. It is not clear to the author whether these problems can be solved by formal means, so any endeavor in this direction may be futile. Accordingly, the reader is advised \emph{not} to read the following too closely. The author feels that he is lured to the same habitual thinking every time he thinks about the problem, and he somehow hopes that the interested reader, with a fresh mind, might be able to come up with a different approach. In fact, such an interested reader might even want to ignore this section (and the one that follows) altogether.
%\erm

In our above argument, the functor that we want to use to relate the derived categories of $\frX^-$ and $\frX^+$ can be written as follows:
\begin{equation}  \label{comp_DA_functor}
	\Dmcoh(\frX^+) \xlaq[\simeq]{(j^+)^*} \cD^-_{{\rm coh},[\geq w],(<w)}(\GrA) \xraq[(?)]{\bD_A} \cD^+_{{\rm coh},[\leq -w],(>-w)}(\GrA)^{\op} \xraq{(j^-)^*}  \Dpcoh(\frX^-)^{\op}
\end{equation}
Here, we have marked the middle arrow with a question mark to indicate that it is not clear whether the functor $\bD_A$ really sends $\cD^-_{{\rm coh},[\geq w],(<w)}(\GrA)$ to $\cD^+_{{\rm coh},[\leq -w],(>-w)}(\GrA)$. Thus, the question mark in \eqref{comp_DA_functor} is a version of \eqref{DA_idea_1} where we restrict to $\Dmcoh$. If it holds, then we see that \eqref{comp_DA_functor} is a composition of three fully faithful functors, and is therefore fully faithful, giving a result along the lines of \eqref{expected_pic_flip_flop}. 

%We do not know whether the last functor $(j^-)^*$ is an equivalence. However, we know that the functor $(j^-)^* : \cD_{[\leq -w],(>-w)}(\GrA) \ra \cD(\frX^+)$ is an equivalence, so its restriction $(j^-)^* : \cD^+_{{\rm coh}[\leq -w],(>-w)}(\GrA) \ra \Dpcoh(\frX^-)$ is fully faithful. Since $A$ is Gorenstein, the functor $\bD_A$ is also fully faithful on $\cD_{{\rm coh}}(\GrA)$. Thus, if one can establish \eqref{DA_idea_1}, or rather the version espressed in (?) here, then the composition \eqref{comp_DA_functor} will be fully faithful.

Whether or not the question mark in \eqref{comp_DA_functor} holds, the functor \eqref{comp_DA_functor} makes sense if we consider
\begin{equation*}  
	\Dmcoh(\frX^+) \xlaq[\simeq]{(j^+)^*} \cD^-_{{\rm coh},[\geq w],(<w)}(\GrA) \xraq{(j^-)^* \circ \bD_A}  \Dpcoh(\frX^-)^{\op}
\end{equation*}
Since duality is local on $\Dmcoh(\frX)$, the last functor may be rewritten as $(j^-)^* \circ \bD_A = \bD_{\frX^-} \circ (j^-)^*$. 
But $\frX^-$ is Gorenstein, so that $\bD_{\frX^-}$ is a contravariant equivalence on $\cD_{{\rm coh}}(\frX^-)$, and we may skip that functor.
In other words, the desired functor is in fact
\begin{equation}  \label{comp_functor_Dmcoh}
	\begin{tikzcd} [column sep = 50]
		\Dmcoh(\frX^+) \ar[r, shift left = 1.5, "\simeq"', "\cL_{[\geq w]} \circ {\bm R}j^+_*"] &  \cD^-_{{\rm coh},[\geq w],(<w)}(\GrA) \ar[l, shift left = 1.5, "(j^+)^*"] \ar[r, "(j^-)^*"] &  \Dmcoh(\frX^-)
	\end{tikzcd}
\end{equation}

Since the duality $\bD_A$ is not involved, we expect that our argument will not need to use the Gorenstein property. In fact, we will show that the duality $\RGIp(A)(a)[1] \cong \bD_Y(\RGIm(A))$ ``almost implies'' that this functor is fully faithful, except that a similar formal problem arise.
To describe the problem, we start with the following straightforward
\blm
Given a full triangulated subcategory $\cE \subset \cD(\GrA)$, then the functor $\CIm :  \cE \ra \cD_{\Imtr}(\GrA)$
is fully faithful if and only if 
\begin{equation*}
	\RHom_A(M, \RGIm(N)) := \RHomcom_A(M, \RGIm(N))_0 \simeq 0 \qquad \text{ for all } M,N \in \cE
\end{equation*}
\elm

We wish to show that the second functor in \eqref{comp_functor_Dmcoh} is fully faithful by verifying this condition.
Thus, let $M \in \cD^-_{{\rm coh},(<w)}(\GrA)$ and $N \in \cD^-_{{\rm coh},[\geq w]}(\GrA)$, we compute
\begin{equation*}  
	\RHomcom_A(M, \RGIm(N)) \stackrel{(?)}{\simeq} \RHomcom_A(M, \RGIm(A)) \otimes_A^{{\bm L}} N \simeq \bD_Y(\RGIp(M))(-a)[-1] \otimes_A^{{\bm L}} N
\end{equation*}
Notice that $\bD_Y(\RGIp(M)) \otimes_A^{{\bm L}} N$ is in $\cD_{>-w} \otimes_A^{{\bm L}} \cD_{[\geq w]} \subset \cD_{>0}$. Therefore, its shift by $a \geq 0$ is zero at weight zero.

We see once again that the argument falls short because of what seems like a formal problem (more precisely, the equality marked with (?) may not hold).
However, notice that our arguments in either the formulations \eqref{DA_idea_1} or \eqref{comp_functor_Dmcoh} work if the $3$-term SOD restricts to $\Dperf(\GrA)$. In the next section, we will see two examples in which it holds.

\brm
One can ``solve'' the problem \eqref{DA_idea_2}(1) by formulating the duality $\bD_A$ in terms of a perfect pairing of DG categories 
\begin{equation*} 
	\cD_{[\geq w]}(\GrA) \times \cD_{[\leq -w]}(\GrA) \raq \cD(k) \, , \qquad 
	(M,N) \mapsto (M \otimes^{{\bm L}}_A N)_0
\end{equation*}
which then restricts to a pairing on the window subcategories $\cD_{[\geq w],(<w)}(\GrA)$ and $\cD_{[\leq -w],(>-w)}(\GrA)$. One can then try to show that it is a perfect pairing.
Unravelling the definition (at least in the way the author formulated it), it corresponds to the functor \eqref{comp_functor_Dmcoh}, except that we do not restrict to $\Dmcoh$.
\erm

\brm
The functor \eqref{comp_functor_Dmcoh} is $\cO_Y$-linear. 
\erm

\brm
By a result \cite{Kaw08} of Kawamata, any two birational projective Calabi-Yau varieties with at most $\bQ$-factorial terminal singularities are connected by a finite number of flops, each of which is obviously a Gorenstein flop, and hence falls into our present setting. 
\erm

\brm
One might try to repeat the arguments in this section for a different version of $\cD(\GrA)$. Some potential candidates are: (1) ${\rm Ind} \, \Dbcoh(\GrA)$, (2) co-derived or contra-derived categories, (3) homotopy category of (unbounded) complexes of finitely generated projectives.
For example, if we work with (3), then problems such as \eqref{bdd_below_complex_weight} seem to not arise (although other problems might arise then).
Accordingly, the theory of relative homological algebra ({\it e.g.} \cite{EJ95}) might be relevant here, although the author is not very familiar with this part of the literature.
If one changes the setting, one should also guarantee that results of Section \ref{flip_flop_sec} can be carried over. 
It is not clear to us whether endeavors in these formal directions will be useful (the author has spent a lot of time trying out various formal modifications, but such effort has not been fruitful so far).
\erm

\section{Some examples}  \label{example_sec}

We will present two examples in which the $3$-term SOD \eqref{three_term_SOD} restricts to $\Dperf(\GrA)$, so that our argument does work. For this to hold, it is certainly necessary that the functor $\cL_{[\geq w]}$ in \eqref{weight_trunc_seq} preserves $\Dperf(\GrA)$.
The following result says that this condition is also sufficient. A proof may be found in \cite{Yeu20b} (a revised version is under preparation, for which the arguments will be simplified).
\bpp
Suppose that the functor $\cL_{[\geq w]}$ in \eqref{weight_trunc_seq} preserves $\Dperf(\GrA)$, then the the $3$-term SOD \eqref{three_term_SOD} restricts to $\Dperf(\GrA)$.
\epp

In this case, the functor \eqref{comp_functor_Dmcoh} restricts to 
\begin{equation}  \label{comp_functor_Dperf}
	\begin{tikzcd} [column sep = 50]
		\Dperf(\frX^+) \ar[r, shift left = 1.5, "\simeq"', "\cL_{[\geq w]} \circ {\bm R}j^+_*"] &  \cD_{\perf,[\geq w],(<w)}(\GrA) \ar[l, shift left = 1.5, "(j^+)^*"] \ar[r, "(j^-)^*"] &  \Dperf(\frX^-)
	\end{tikzcd}
\end{equation}
Moreover, the arguments in the previous section then guarantees that, if the duality $\bD_Y(\RGIm(A)) \cong \RGIp(A)(a)[1]$ holds for some $a \geq 0$, then \eqref{comp_functor_Dperf} is fully faithful.

In general, it seems quite difficult to describe $\cL_{[\geq w]}$ explicitly, as it is defined in terms of the Neeman-Brown representability theorem. The closest to an explicit formula that we managed to get is the following tensor form
\begin{equation}  \label{cL_tensor_form}
	\cL_{[\geq w]}(M) \, = \, M_{\geq w} \otimes^{{\bm L}}_{\cF_{\geq w}} \cF
\end{equation}
Here, $\cF$ is the small pre-additive ({\it i.e.,} $\Ab$-enriched) category given by $\Ob(\cF) = \bZ$ and $\cF(i,j) = A_{i-j}$. Then a graded $A$-module is the same as a right module over $\cF$ (recall that a right module over a small pre-additive category $\cC$ means an additive functor $\cC^{\op}\ra \Ab$).  Let $\cF_{\geq w} \subset \cF$ be the full subcategory on the object set $\bZ_{\geq w} \subset \bZ$. Then the inclusion functor $\cF_{\geq w} \ra \cF$ gives rise to restriction and induction functors between the module categories. In particular, for any $M \in \Mod(\cF)$, we may restrict it to $\cF_{\geq w}$ to obtain $M_{\geq w} \in \Mod(\cF_{\geq w})$, and then tensor it back to $\cF$ to obtain $M_{\geq w} \otimes_{\cF_{\geq w}} \cF \in \Mod(\cF)$. The right hand side of \eqref{cL_tensor_form} is its derived functor. %Indeed, \eqref{cL_tensor_form} is simply an unravelling of (the DG-categorical version of) the Neeman-Brown representability theorem.

The description \eqref{cL_tensor_form} seems to be of little use in the two examples that we present below. In both cases, the functor $\cL_{[\geq w]}$ is described by some ad-hoc arguments specific to the form of the graded ring $A$ in question.

The first example is the polynomial ring
\begin{equation}  \label{poly_ring}
	A = k[x_1,\ldots,x_p,y_1,\ldots,y_q] \, , \qquad \deg(x_i) > 0, \, \, \deg(y_j)<0
\end{equation}
Write $A = A^+ \otimes A^-$, where $A^+ = k[x_1,\ldots, x_p]$  and $A^- = k[y_1,\ldots, y_q]$. 
Then we have
\begin{equation*}
	\cD_{[\geq w]}(\Gr(A^+)) \otimes A^- \subset \cD_{[\geq w]}(\Gr(A))
	\qquad \text{and} \qquad 
	\cD_{< w}(\Gr(A^+)) \otimes A^- \subset \cD_{< w}(\Gr(A))
\end{equation*}
This allows us to compute the weight truncation of objects of the form $M = M^+ \otimes A^-$ for $M^+ \in \cD(\Gr(A^+))$. Namely, we have
\begin{equation*}
	\cL_{[\geq w]}^A( M^+ \otimes A^- ) \, = \, \cL_{[\geq w]}^{A^+}( M^+) \otimes A^- 
\end{equation*}

To show that $\cL_{[\geq w]}^A$ preserves $\Dperf(\GrA)$, it suffices to verify it on the split-generating objects $\{A(-i)\}_{i \in \bZ}$, each of which is of the form $M^+ \otimes A^-$ for $M^+ = A^+(-i)$. Hence it suffices to compute $\cL_{[\geq w]}^{A^+}( A^+(-i) )$. 
Since $A^+$ is non-negatively graded, we in fact have $\cD_{[\geq w]}(\Gr(A^+)) = \cD_{\geq w}(\Gr(A^+))$, so that the weight truncation is easy to describe:
\blm  \label{wtr_poly_lem}
For any $M \in \Gr(A^+)$, we have
\begin{equation*}
	\cL_{[\geq w]}^{A^+}( M ) = M_{\geq w} \, , \qquad \text{and} \qquad 
	\cL_{< w}^{A^+}( M ) = M / M_{\geq w}
\end{equation*}
\elm
As a result, we have
\begin{equation*}
	\cL_{[\geq w]}^{A}( A(-i) ) = ((A^+)_{\geq w-i} \otimes A^-)(-i)
\end{equation*}
Since $A^+$ is smooth, we have $(A^+)_{\geq w-i} \in \Dperf(\Gr(A^+))$. This shows the following
\bpp  \label{poly_perfect_wtr}
For the polynomial ring \eqref{poly_ring}, weight truncation $\cL_{[\geq w]}$ preserves $\Dperf(\GrA)$.
\epp

\brm
Our argument here works whenever a graded ring can be written as $A = A^+ \otimes A^-$ where $A^+$ is non-negatively graded and $A^-$ is non-positively graded. In this case, if $A^+$ is smooth, then $\cL_{[\geq w]}$ preserves $\Dperf(\GrA)$. Interestingly, if $A^-$ is smooth, then by the arguments of \cite{HL15}, one can show that $\cL_{[\geq w]}$ preserves $\Dbcoh(\GrA)$.
\erm

Proposition \ref{poly_perfect_wtr} shows that our approach in the previous section works perfectly for the polynomial ring \eqref{poly_ring} (notice that we have established the duality $\bD_Y(\RGIm(A)) \cong \RGIp(A)(a)[1]$ directly in Remark \ref{poly_ring_RGam_rmk} without applying the main Theorems of Section \ref{flip_flop_sec}). 
We will now compute the corresponding functor \eqref{comp_functor_Dperf}. For simplicity, we will now assume $w = 0$.  
Let $\eta^+ = \sum_{i=1}^p \deg(x_i)$ and $\eta^- = -\sum_{j=1}^p \deg(y_j)$, then $a = \eta^- - \eta^+$.

By the computation in Remark \ref{poly_ring_RGam_rmk}, we see that
\begin{equation*}
	\RGIp(A) \in \cD_{\leq -\eta^+}(\GrA) \qquad \text{and} \qquad \RGIm(A) \in \cD_{\geq \eta^-}(\GrA) 
\end{equation*}
Thus, in particular, for $0 \leq i < \eta^+$, the objects $A(-i)$ are in the window subcategory $\cD_{[\geq 0],(<0)}(\GrA)$. 
This implies that 
\begin{equation}  \label{functor_poly_range_1}
	\parbox{40em}{For $0 \leq i < \eta^+$, the functor \eqref{comp_functor_Dperf} for $w=0$ sends $\cO_{\frX^{+}}(-i)$ to $\cO_{\frX^{-}}(-i)$.}
\end{equation}

In fact, this uniquely characterize the functor \eqref{comp_functor_Dperf}.
Namely, the arguments of \cite{HL15} can be employed to give the following characterization of the window subcategory by the ``window length'' $\eta^+$, which then implies that $\{\cO_{\frX^{+}}(-i) \}_{0 \leq i < \eta^+}$ split generates $\Dperf(\frX^+)$.
\bpp  \label{poly_window_length}
For the polynomial ring \eqref{poly_ring}, we have $\cD_{\perf,(<w)}(\GrA) = \cD_{\perf,[<w+\eta^+]}(\GrA)$. As a result, the subcategory $\cD_{\perf,[\geq 0],(<0)}(\GrA)$ is split generated by $\{A,A(-1),\ldots, A(-\eta^+ +1) \}$. 

Similarly, we have $\cD_{\perf,(>-w)}(\GrA) = \cD_{\perf,[>-w-\eta^-]}(\GrA)$, so that $\cD_{\perf,[\leq 0],(>0)}(\GrA)$ is split generated by $\{A,A(1),\ldots, A(\eta^- -1) \}$. 
\epp

By our above description of $\cL_{[\geq 0]}$ in Lemma \ref{wtr_poly_lem}, we also have
\begin{equation}  \label{functor_poly_range_2}
	\parbox{40em}{For $i > 0$, the functor \eqref{comp_functor_Dperf} for $w=0$ sends $\cO_{\frX^{+}}(i)$ to the associated sheaf on $\frX^-$ of the graded module $((A^+)_{\geq i} \otimes A^-)(i)$}
\end{equation}

The image of $\cO_{\frX^{+}}(i)$ for $i \leq -\eta^+$ seems to be more difficult to explicate, since $\cL_{[\geq 0]}(\RGIp(A(i)))$ does not vanish in this case, so the answer invovles a certain cone. We will skip the details.

\brm
The case $\deg(x_i)=1$ and $\deg(y_j)=-1$ corresponds to a standard flip/flop (we then have $\frX^{\pm} = X^{\pm}$ and $\eta^+ = p$, $\eta^- = q$). In this case, the functor \eqref{comp_functor_Dperf} is isomorphic the Fourier-Mukai functor ${\bm R}q_* \circ {\bm L}p^*$ with respect to the standard span $X^- \xla{q} \widetilde{X} \xra{p} X^+$. Indeed, one can simply verify this on the split generators $\{\cO_{X^{+}}(-i) \}_{0 \leq i < p}$ of $\Dperf(X^+)$. The author thanks Alexander Kuznetsov for help with this remark.
\erm

We now move to our second class of examples, which are the graded rings associated to a class of $3$-fold flips of type A, worked out by Brown and Reid (see \cite[Section 11.2]{Rei}). 

Fix positive integers $d,e,\alpha,\beta,\lambda,\mu$ with ${\rm gcd}(\lambda,\mu)=1$. Take%
\footnote{We have renamed the generators of \cite[Section 11.2]{Rei} according to $(t,u,x_0,x_1,y_0,y_1)\ra (y_3,z,y_1,x_1,y_2,x_2)$, in order to conform with our standing convention that $\deg(x_i)>0$, $\deg(y_i)<0$ and $\deg(z_i)=0$.}
%
%\begin{equation} 
%	A = k[t,u,x_0,x_1,y_0,y_1]/(f_1,f_2)
%\end{equation}
%
\begin{equation}  \label{BR_example}
	A = k[x_1,x_2,y_0,y_1,y_2,z]/(f_1,f_2)
\end{equation}
with degrees and relations given by
\begin{gather*} 
	 \deg(x_1) = \lambda \, , \quad \deg(x_2) = \mu \, , \quad \deg(y_1) = -\mu \, , \quad \deg(y_2) = -\lambda - \mu e \, , \quad \deg(y_3) = -1 \, , \quad \deg(z) = 0 \\
	f_1 = x_1y_2 - y_1^e z^{\alpha} - y_3^{\mu e} \, , \quad
	f_2 = y_1 x_2 - z^{\beta} - x_1^d y_3^{\lambda d}
\end{gather*}

This is an example of a complete intersection, so we will work in that more general context below.
\begin{equation}  \label{poly_pres}
	A = k[x_1,\ldots,x_p,y_1,\ldots,y_q,z_1,\ldots,z_r]/(f_1,\ldots,f_s) \, , \qquad \deg(x_i) > 0, \, \, \deg(y_i)<0 , \, \, \deg(z_i)=0
\end{equation}
Our method will work under the following assumptions (which is satisfied by the example \eqref{BR_example} of Brown and Reid):
\begin{equation}  \label{CI_assump_1}
	\parbox{40em}{$A$ is a complete intersection ({\it i.e.,} $\dim(A) = p+q+r-s$) and $\deg(f_i) \leq 0$ for each $1 \leq i \leq s$.}
\end{equation}

If we write $C = k[x_1,\ldots,x_p,y_1,\ldots,y_q,z_1,\ldots,z_r]$, then this condition guarantees that the Koszul resolution $K_{\bullet}(C,f_1,\ldots,f_s) \simeq A$ gives a resolution of $A$ by free graded $C$-modules, all of which are generated in non-positive degrees. Thus, it gives a non-positive presentation of $A$, in the sense of the following
\bdf
A Noetherian graded algebra $A$ is said to \emph{have a non-positive presentation} if there exists a map $C \ra A$ from a polynomial graded algebra $C$ such that, as an object in $\cD(\Gr(C))$, we have $A \in \cD_{[\leq 0]}(\Gr(C))$.
\edf

\bpp
If $A$ has a non-positive presentation, then $\cL_{[\geq w]}$ preserves $\Dperf(\GrA)$.
\epp

\bpf
In general, we have $\cD_{[\geq w]}(\Gr(C)) \otimes^{{\bm L}}_C A \subset \cD_{[\geq w]}(\GrA)$. The condition on non-positive presentation also guarantees that $\cD_{<w}(\Gr(C)) \otimes^{{\bm L}}_C A \subset \cD_{<w}(\GrA)$. Hence we have $\cL_{[\geq w]}^A(M\otimes^{{\bm L}}_C A) \cong \cL_{[\geq w]}^C(M) \otimes^{{\bm L}}_C A$ for all $M \in \cD(\Gr(C))$. Apply this to $M = C(-i)$ for $i \in \bZ$ to conclude the proof.
\epf

\brm
With the help of computer programs, given any finite map $C \ra A$ from a polynomial graded algebra, one should be able to compute a free resolution of $A$ as a graded $C$-module, from which one should be able to see if the given map $C \ra A$ gives a non-positive presentation. In fact, it suffices to compute the weights of $\Tor^C_{\bullet}(A,k[z_1,\ldots,z_r])$ to verify that $C \ra A$ gives a non-positive presentation.
\erm

Combining with the previous results, we then have the following
\bthm \label{BR_ff}
For the class of 3-fold flips described by \eqref{BR_example}, there is a fully faithful functor $\Dperf(\frX^+) \ra \Dperf(\frX^-)$ given by \eqref{comp_functor_Dperf}.
\ethm

\brm
Let $X$ be a $\bQ$-Gorenstein variety such that $\cO_X(iK_X)$ is maximal Cohen-Macaulay for each $i \in \bZ$ (this is satisfied, {\it e.g.}, if $X$ has at most log terminal singularities, see \cite[Corollary 5.25]{KM98}), one may define its \emph{Gorenstein root} to be the Deligne-Mumford stack $\frX := [ \, \underline{\Spec}_X(\bigoplus_{i \in \bZ} \cO_X(iK_X)) \, / \, \bG_m \, ]$. Theorem \ref{BR_ff} then relates the derived categories of the Gorenstein roots of the two sides of the flip.
\erm

%\brm
%In the example of \eqref{BR_example}, we have seen that $\cL_{[\geq w]}$ preserves $\Dperf(\GrA)$. However, it seems that $\cL_{[\leq -w]}$ does not preserve $\Dperf(\GrA)$. Instead, by the arguments of \cite{HL15}, it preserves $\Dbcoh(\GrA)$. Thus, we have a window $\cD_{\perf,[\geq 0],(<0)}(\GrA)$ for $\Dperf(\frX^+)$ in the positive direction, and a window $\cD^b_{{\rm coh},[\leq 0],(>0)}(\GrA)$ for $\Dbcoh(\frX^-)$ in the negative direction. Since the structure sheaf of the unstable locus in the negative direction is not perfect in the ambient stack $\frX$, there is apparently no ``length characterization'' of the window  $\cD^b_{{\rm coh},[\leq 0],(>0)}(\GrA)$ for $\Dbcoh(\frX^-)$. However, its intersection with $\Dperf(\GrA)$, {\it i.e.,} $\cD_{\perf,[\leq 0],(>0)}(\GrA)$, although not a window for $\Dperf(\frX^-)$, does have a length characterization. Thus, it seems that one can still apply the arguments in \cite{HL15} to establish Theorem \ref{BR_ff}. 
%\erm

One can also compute the functor \eqref{comp_functor_Dperf} in the present setting. We will work with \eqref{poly_pres} satisfying \eqref{CI_assump_1}. Let $C = k[x_1,\ldots,x_p,y_1,\ldots,y_q,z_1,\ldots,z_r]$, and let $C^+ = k[x_1,\ldots,x_p]$ and $C^- = k[y_1,\ldots,y_q,z_1,\ldots,z_r]$, so that $C = C^+ \otimes C^-$. Let $\eta^+ = \sum_{i=1}^p \deg(x_i)$. Then we have
\begin{equation}  \label{functor_CI_range_1}
	\parbox{40em}{For $0 \leq i < \eta^+$, the functor \eqref{comp_functor_Dperf} for $w=0$ sends $\cO_{\frX^{+}}(-i)$ to $\cO_{\frX^{-}}(-i)$.}
\end{equation}
\begin{equation}  \label{functor_CI_range_2}
	\parbox{40em}{For $i > 0$, the functor \eqref{comp_functor_Dperf} for $w=0$ sends $\cO_{\frX^{+}}(i)$ to the associated sheaf on $\frX^-$ of the complex of graded modules $((C^+)_{\geq i} \otimes C^-) \otimes^{{\bm L}}_C A(i)$}
\end{equation}

We may also want an analogue of Proposition \ref{poly_window_length}. For that, we will need to replace the assumption \eqref{CI_assump_1} with the following stronger one (which is satisfied by \eqref{BR_example}):
\begin{equation}  \label{CI_assump_2}
\parbox{40em}{$\dim(A) = p+q+r-s$, $\dim(A/I^+) = q+r-s$, and $\deg(f_i) \leq 0$ for each $1 \leq i \leq s$.}	
\end{equation}
Then we have the following

\bpp
For $A$ in \eqref{poly_pres} satisfying \eqref{CI_assump_2}, we have $\cD_{\perf,(<w)}(\GrA) = \cD_{\perf,[<w+\eta^+]}(\GrA)$. As a result, the subcategory $\cD_{\perf,[\geq 0],(<0)}(\GrA)$ is split generated by $\{A,A(-1),\ldots, A(-\eta^+ +1) \}$. 
\epp

\textbf{Acknowledgement.} The author thanks Daniel Halpern-Leistner and Valery Lunts for helpful discussions. He thanks Alexander Kuznetsov for careful reading and helpful comments and suggestions. He thanks Chen Jiang for help with Theorem \ref{KV_vanishing}.


\begin{thebibliography}{9}
%----A----




%----B----
\bibitem{BBD82}
A. Beĭlinson, J. Bernstein, and P. Deligne,
\textit{Faisceaux pervers}, 
Analysis and topology on singular spaces, I (Luminy, 1981), 5--171, 
Astérisque, 100, Soc. Math. France, Paris, 1982. 

\bibitem{BH93}
W. Bruns, and J. Herzog, 
\textit{Cohen-Macaulay rings}, 
Cambridge Studies in Advanced Mathematics, \textbf{39}. Cambridge University Press, Cambridge, 1993.

\bibitem{BFK19}
M. Ballard, D. Favero, and L. Katzarkov, 
\textit{Variation of geometric invariant theory quotients and derived categories}, J. Reine Angew. Math. \textbf{746} (2019), 235--303.
%----C----



%----D----




%----E----
\bibitem{EJ95}
E.E. Enochs, and O.M.G. Jenda, \textit{Gorenstein injective and projective modules}, Math. Z.
\textbf{220} (1995), 611--633

%----F----



%----G----




%----H----
\bibitem{Har77}
R. Hartshorne, 
\textit{Algebraic Geometry}, Graduate Texts in Mathematics, \textbf{52}. Springer-Verlag, New York-Heidelberg, 1977.

\bibitem{HL15}
D. Halpern-Leistner, 
\textit{The derived category of a GIT quotient},
J. Amer. Math. Soc. \textbf{28} (2015), 871--912. 


%----I----



%----J----

%----K----
\bibitem{KM98}
J. Koll\'{a}r, and S. Mori, \textit{Birational geometry of algebraic varieties}, With the collaboration of C. H. Clemens and A. Corti. Translated from the 1998 Japanese original. Cambridge Tracts in Mathematics, 134. Cambridge University Press, Cambridge, 1998.

\bibitem{Kaw08}
Y. Kawamata, 
\textit{Flops connect minimal models},
Publ. Res. Inst. Math. Sci. \textbf{44} (2008), 419--423.

\bibitem{KKM87}
Y. Kawamata, K. Matsuda, and K. Matsuki, \textit{Introduction to the minimal model problem}, Algebraic geometry (Sendai, 1985) (Advanced Studies in Pure Mathematics) Volume 10, North-Holland, 1987, 283--360
%----L----


%----M----
\bibitem{May01}
J. P. May, 
\textit{The additivity of traces in triangulated categories},
Adv. Math. \textbf{163} (2001), no. 1, 34--73.

%----N----




%----O----
\bibitem{Orl09}
D. Orlov, \textit{Derived categories of coherent sheaves and triangulated categories of singularities}, Algebra, arithmetic, and geometry: in honor of Yu. I. Manin. Vol. II, 503--531, Progr. Math., \textbf{270}, Birkhäuser Boston, Inc., Boston, MA, 2009.


%----P----


%----Q----



%----R----
\bibitem{Rei}
M. Reid, \textit{Graded rings and birational geometry},
notes available online.



%----S----




%----T----



%----U----



%----V----


%----W----



%----X----



%----Y----
\bibitem{Yeu20a}
W.K. Yeung, 
\textit{Homological flips and homological flops}, {\tt arXiv:1907.06190}.

\bibitem{Yeu20b}
W.K. Yeung,
\textit{Weight truncation for wall-crossings in birational cobordisms}, {\tt arXiv:2001.10431}

\bibitem{Yeu20c}
W.K. Yeung,
\textit{Grothendieck duality and Greenlees-May duality on graded rings}, {\tt arXiv:2001.08795}.


%----Z----














\end{thebibliography}
\end{document}